\documentclass[12pt]{amsart}
\usepackage{amsmath,amssymb}
\title[The Uniqueness Of Integrals ]{The Uniqueness Of Integrals On Hopf Algebras\\
 A categorical approach} 

\author{Ph\`ung H$\grave{\mbox{\^o}}$ Hai }
\address[PHH]{Hanoi Institute of Mathematics, P.O Box 631, 10000 Boho, Hanoi, Vietnam}
\author{Nguy$\tilde{\mbox{\^e}}$n Huy Hung}
\address[NHH]{ Hanoi Pedagogical University II, Melinh, Vinhphuc, Vietnam }

\newtheorem{thm}{Theorem}
\newtheorem{lemma}[thm]{Lemma}
\newtheorem{corollary}[thm]{Corollary}
\newtheorem{pro}[thm]{Proposition}

\def\eee{\rule{.75ex}{1.5ex}\\[0ex]}
\begin{document}

\maketitle
\bibliographystyle{plain}

\def\phi{\varphi}
\def\lhom{\mbox{\rm\sffamily lhom}}
\def\rhom{\mbox{\rm\sffamily rhom}}
\def\ComodH{\mbox{\rm\sffamily Comod-$H$}}
\def\HMod{\mbox{\rm\sffamily $H^*$-Mod}}
\def\ModA{\mbox{\rm\sffamily Mod-$A$}}
\def\Mod{\mbox{\rm\sffamily Mod}}
\def\db{\mbox{\rm\sffamily db}}
\def\ev{\mbox{\rm\sffamily ev}}
\def\id{\mbox{\rm\sffamily id}}

\def\dim{\mbox{\rm dim}}
\def\Hom{\mbox{\rm Hom}}
\def\End{\mbox{\rm End}}
\def\ot{\otimes }
\def\loma{\longrightarrow}
\def\rint{{\textstyle\int_r}}
\def\lint{{\textstyle\int_l}}

\begin{abstract}
  We propose a new method to 
investigate the dimension of the space of integrals on a Hopf algebra
 $H$ and other properties of $H$-comodules. 
\end{abstract}

\section*{Introduction}
The study of integrals is an important part in the theory of Hopf algebras.
If the integral does not vanish at the unit, the Hopf algebra is cosemisimple.
In general, if the space of integrals doen not vanish, many interesting properties
of the Hopf algebra can be derived, for instance, the bijectivity of the antipode,
the finiteness of injective envelopes, etc...

The study of the space of integrals on a Hopf algebras was initiated by Sweedler,
being motivated by the uniquenness of the Haar integral on locally compact groups.
Sweedler proved the existence and 
uniqueness up to a constant 
of a non-zero integral on any finite dimensional Hopf algebra and asked whether this
remains true in infinit dimensional case. This question was answered by Sullivan \cite{sul} in 1971.
He showed that the space of integrals has dimension not preceeding 1, in other words,
if a non-zero integral exists then it is uniquely determined upto a constant.
Since then there have been several other proofs of this result \cite{stefan1,bdgn}.

The aim of this paper is to give a new proof of the uniqueness of integrals on 
Hopf algebras and properties 
of right $H$-comodules. Our method bases on a theorem of 
Gabriel-Popescu characterizing 
Grothendieck categories (see, e.g., \cite{stenstrom}). We started with the known fact that a Hopf algebra
over a field is completely determined by the category of its, say, right comodules.
The question is then whether can we express the properties of the integral space
and related properties in terms of the comodule category. It turns out that a
 non-zero integral on a Hopf algebra $H$ exists if an only if $H$ is a generator in the category
Comod-$H$ of its right comodules. This category is a Grothendieck category and hence,
by the theorem of Gabriel-Popescu, is the quotient of a module category. The latter is
indeed the category of (left) $H^*$-modules.

The paper is organized as follows. After recalling some basic fact about Hopf algebras,
their comodules and their duals we mention some isomorphisms which play essential
role in our work. Then we recall the notion of Grothendieck categories and the
theorem of Gabriel-Popescu characterizing these categories (Section 1). We prove
in Section 2 the uniqueness of integral. We also mention some interesting consequences
of our methods: Corollaries \ref{inj} and \ref{lem_gamma}.
In Section 3 we describe the structure of the functor $U$ and derive from this the bijectivity
of the antipode. Finally we prove that the injective envelope of $k$ is finite dimensional.

Although our main results are known, our method seems to be new and we hope that
this method can be generalized to the case of Hopf algebras (or algebroids) over a ring.

\section{Preliminaries}
\subsection{Integrals}
Let $k$ be a fixed field and $H$ be a Hopf algebra defined over $k$.
A right integral on $H$ is a linear functional $\chi:H\longrightarrow k$, subject to the following
equation
\begin{equation}\sum_{(a)}\label{left_int} \chi(a_{(1)})a_{(2)}=\chi(a),\quad \forall a\in H.
\end{equation}
The space of right integrals on $H$ is denoted by $\rint$.
Consider $k$ as the trivial $H$-comodule by means of the unit map, we see that
$\rint\cong\Hom^H(H,k)$. Analogously, a left integral on $H$ is a linear map
$\varphi:H\longrightarrow k$, satisfying the  equation 
$\sum_{(a)} a_{(1)}\chi(a_{(2)})=\chi(a),\quad \forall a\in H.$ The space of left integrals is
denoted by $\lint$, consider $k$ as a left $H$-comodule, we have $\lint\cong{}^H\Hom(H,k)$.

\subsection{The dual algebra $H^*$}
Let  $M$ be a right $H$-comodule with coaction $\rho$ and $X$ be a $k$-space. 
The coproduct on $H$ induces a coaction on $X\ot H$, trivial on $X$:
$x\ot h\loma \sum_{(h)}x\ot h_{(1)}\ot h_{(2)}.$ We denote this comodule by 
$(X)\ot H$ to emphasize that the coation does not involve $H$. The map 
\begin{equation}
 \Hom ^H (M, (X) \otimes _k H) \longrightarrow \Hom _k (M,X), f \mapsto (id \otimes \epsilon )
\cdot f  \label{eq0}
\end{equation} 
is an isomorphism with the inverse map
given by  $h \mapsto (h \otimes id) \cdot \rho$. 

In particular, for $X = k$ we have an isomorphism $\Hom ^H (M,H) \cong \Hom _k (M,k)=:M^*$.
For $M=H$ we have an isomorphism $\Hom ^H (H,H) \cong H^*$ which is in fact an algebra 
anti-isomorphism, with respect to the convolution product on $H^*$ is given by 
$$(f*g)(h):=\sum_{(h)} f(h_{(1)}) g(h_{(2)}).$$

\subsection{Rational comodules}
For any $N$ in $\ComodH $ with the coaction map $\rho : N \longrightarrow N \otimes 
H, \rho (n)= \sum_{(n)} n_{(1)} \otimes n_{(2)} $. Contruct an  $H^*$-module structure on $N$ 
as follows: 
$$\forall \xi \in H^*, \forall n \in N,  \xi * n:= \sum_{(n)} n_{(1)} \xi (n_{(2)})$$
In this way every right $H$-comodule has a structure left $H^*$-module. 
This contruction indeed defines a functor from $\ComodH$ to $H^*$-\Mod\ which 
is fully faithful and exact. An $H^*$-module which is isomorphic to a module obtained
in this way is called {\em rational} module. It is known that submodules and quotient 
modules of a rational module are rational, and sum of rational submodules is rational.

Not every $H^*$-module is rational. However, we have 
a left exact functor $Rat$ from $H^*$-\Mod\ to $\ComodH $  \cite{sweedler1}, assigning to any 
left $H^*$-module $M$ the largest rational submodule of $M$ in $\ComodH$.

\subsection{Tensor product and internal homs}
The category of right $H$-comodules is a monoidal category  
with the usual tensor product over $k$
 and the unit object is $k$. An object $ X \in \ComodH $ 
defines functors $X \otimes -$ and $- \otimes X$ 
from $\ComodH $ to itself. The right adjoint functors to these 
functors are called the left and the right internal hom 
functors of $X$ and denoted by 
$\rhom  (X,-)$ and $\lhom    (X,-)$, respectively. Thus we have natural isomorphisms 
$$\Hom^H(X \otimes Y, Z) \cong \Hom^H(Y, \rhom (X,Z))$$
$$\Hom^H(Y \otimes X, Z) \cong \Hom^H(Y, \lhom  (X,Z))$$ 
for all $Y, Z $ in $\ComodH $. 

Let $N$ be a finite dimensional right $H$-comodule then $N^* := \Hom_k (N,k)$
 is also a right comodule with 
the coaction given by the equation
$$\rho (\phi )(x) := \phi _{(0)} (x) \phi _{(1)} = \phi (x_{(0)}) S(x_{(1)}), x \in N, \phi \in N^* . $$
The map ev$:M^* \otimes M \rightarrow k, \phi  \otimes x \mapsto \phi (x) 
$ is a morphism of $H$-comodules. The pair 
$(M^*$ , ev) is called a left dual to $M$; it is defined uniquely up to isomorphisms.
 There exists a map
$\db:k \rightarrow M \otimes M^* $, defined by the conditions
 $(\ev \otimes \id_{M^*})(\id_{M^*} \otimes \db) = \id_{M^*}$ and 
$(\id_M \otimes \ev)(\db \otimes \id_M) =  \id_{M}$, which is also a 
comodule morphism. The notion of right dual is defined analogously;
 for instance, $(M, \ev)$ is the right dual to $M^* .$ 

	We see that the left dual to a finite dimensional comodule always exists.
 If the antipode is bijective then the right dual to any 
finite dimensional comodule also exists. The converse is also true, if the right dual
to any finite dimensional comodule exists then the antipode is bijective, see, e.g.,
\cite{schauen}. We shall need the following isomorphisms, 
given by manipulating the morphism $\ev$ and $\db$:
 for any finite dimensional comodule $N$, 
\begin{equation}\Hom^H (M \otimes N, P) \cong \Hom^H (M, P \otimes N^*) \label{eq14} 
\end{equation}
\begin{equation}\Hom^H (M, N \otimes P ) \cong \Hom^H (N^* \otimes M, P)  \label {eq15}
\end{equation} 
i.e., $\lhom  (N, P) \cong P \otimes N^*$ (as $H$-comodules). 
 
\subsection{Doi's and Sweedler's isomorphisms} 
In the previous subsection, except for the notion of dual comodules we haven't
use the antipode. In fact, the distinguished r\^ole of the antipode can be expressed
by the isomophisms below.
 For any right $H$-comodule $V$, Yu. Doi discovered the following 
isomorphism of $H$-comodules
 \begin{equation}V \otimes H \cong (V) \otimes H \label{eq2} 
\end{equation}
via the map $v \otimes h \mapsto  \sum v_{(0)} \otimes v_{(1)}h $, and its inverse is given by
$v \otimes h \mapsto  \sum v_{(0)} \otimes S(v_{(1)}) h, $ \cite[1.3 Cor1]{doi1}.

Another isomorphism was established by Sweedler \cite[5.1.3]{sweedler1}:
\begin{equation}\lint \otimes H \cong Rat(H^*). \label{eq5}\end{equation} 
given by $\phi \ot h \mapsto \phi_h: \phi_h(a)=\phi(aS(h)).$ A similar isomorphism
exists for right integrals.

Given an element $\phi\in\lint$, define the
 map $\phi^*:H\longrightarrow H^*$, $\phi^*(h)(a)=\phi(hS(a))$. This  is a morphism
of right $H^*$-modules. Indeed, for all $\xi\in H^*$:
$$\phi^*(h*\xi)(a)=\sum_{(a)}\xi(h_{(1)})\phi(h_{(2)}S(a))=\sum_{(a)}\phi(hS(a_{(1)})\xi(a_{(2)})
=(\phi^*(h)*\xi)(a).$$
Thus, if $\lint\neq 0$, we can chose a left integral $\phi$ such that $\phi^*$ is non-zero.
Consequently $H^*$ considered as right module on itself contains non-zero rational
submodule. Using Sweedler's isomorphism above (for right integrals), we conclude
that $\rint\neq 0$. We therefore have proved (bearing in mind the symmetry
between the notions of left and right integrals)
\begin{lemma}\label {thmbs1} Let $H$ be a Hopf algebra over a field $k$. Then $\rint\neq 0$
if and only if $\lint\neq 0$.
\end{lemma}
\subsection{$H$-Comod is a Grothendieck category}
By definition, a Grothendieck category is an abelian category in which direct limits exist and
preserve left exact sequence and there exists a generator. An example of Grothendieck category
is the category of modules over a ring. In fact, a theorem of Gabriel-Popescu 
\cite[4.1]{stenstrom} states that
a Grothendieck category is the quotient of a module category. More explicitly, the following is 
known. Let $C$ be a generator in the Grothendieck category $\mathcal  C$. Let $A:=\End(C)$. The
functor $T=\Hom(C,-):\mathcal  C\longrightarrow \mbox{Mod-A}$ is then left exact and
fully faithful. There 
exists a left adjoint to this functor, say $U$, which is exact and by means of which $\mathcal  C$ is a
quotient category of Mod-$A$. The adjointness is
\begin{equation}\label{eq_adj} \Hom_{\mathcal  C}(U(M), X)\cong \Hom_A(M,\Hom(C,X)).
\end{equation}
where $M$ is an $A$-module and $X\in\mathcal  C$. It follows from the faithfulness of $\Hom(C,-)$
that 
\begin{equation} U(\Hom(C,X))\cong X, \quad \forall X\in\mathcal C.\end{equation}
\section{The uniqueness of the integral}
Our method of proving the uniqueness of integrals is first to show hat the
category of $H$-comodules is a Grothendieck category with $H$ being a
generator, then to apply the theorem of Gabriel-Popescu to construct an
exact functor $U$ from the category of $H^*$-modules to the category of 
$H$-comodules.
\begin{lemma} \label{lem3}
Let $H$ be a Hopf algebra over a field $k$ with a non-zero right integral. Then $H$ is a 
generator in $ \ComodH $.
\end{lemma}
{\it Proof}. We have according to (\ref{eq15}) and (\ref{eq2})
\begin{eqnarray*} \Hom ^H (H, M )& \cong & \Hom ^H (M^* \otimes H, k) \\
& \cong & \Hom ^H( (M^*) \otimes H, k)  \\
& \cong & M^{**} \otimes  \Hom ^H (H, k) \\
&\neq &0\end{eqnarray*}Thus  $H$ is a generator in $\ComodH $.\eee

The category  $\ComodH $ is a cocomplete
 abelian category, in which direct limits are exact. According to Lemma 
\ref{lem3}, $\ComodH $ is a Grothendieck category 
with $H$ being a generator. Putting $A:=\Hom ^H (H,H)$, we have functor
\begin{equation} T: \ComodH  \longrightarrow \ModA, X \mapsto 
\Hom ^H (H,X).  \label{eq6}\end{equation}
 $T$ is full and faithful and according to Gabriel-Popescu's Theorem  for $\ComodH $, it has a left 
adjoint functor $U: \ModA \longrightarrow \ComodH, $ which is exact:
$$\Hom^H(U(M),N)\cong \Hom_A(M,T(N)).$$
Remember that  $A$ and $H^*$ are anti-isomorphic, so every right $A$-module 
is a left  $H^*$-module and conversely. The explicit action of $H^*$ on 
$\Hom ^H (H,N)$ is given by 
\begin{equation} \xi *f (h):= \sum_{(h)} \xi (h_{(1)}) f(h_{(2)}), 
\mbox{  for }\xi \in H^*, f \in \Hom ^H (H,N), h \in H \label{eq8}\end{equation}  
Thus we have 
\begin{equation}\Hom ^H (U(M), N) \cong \Hom _{H^*} (M,\Hom(H,N)). \label{eq9} \end{equation}
for any  left $H^*$-module $M$ and right $H$-comodule $N$.   
In particular, for any $N$ in $\ComodH $, 
\begin{equation}U(\Hom ^H (H,N)) \cong N. \label{eq7}\end{equation}

We have proved
\begin{pro} \label{pro4}
Let $H$ be a Hopf algebra with a non-zero right integral  and 
$$T: \ComodH  \longrightarrow \HMod, X \mapsto \Hom ^H (H,X)$$
be a functor from $\ComodH $ to  $\HMod$. Then there exists a functor 
$U: \HMod \longrightarrow \ComodH $, which is a left adjoint functor of 
$T$ and $U$ is exact.  \eee
\end{pro}
Since $H$ is injective in Comod-$H$, using standard arguments we deduce
\begin{corollary}\label{inj} With the assumption of Proposition \ref{pro4}, $H^*$ is
injective in the category $H^*$-\Mod.\end{corollary}

On the other hand, setting $N = H$ in (\ref{eq9}), we have
\begin{eqnarray*} \Hom ^H (U(M), H) & \cong  &\Hom _{H^{*}} (M, \Hom ^H (H,H)) \\
& \cong &  \Hom _{H^{*}} (M, H^{*}),\end{eqnarray*}
that is, for any M in $\ComodH$,
\begin{equation}
\Hom ^H (U(M), H) \cong \Hom _{H^{*}} (M, H^{*}).   \label{eq10}
\end{equation}  
For  $M$ in $\ComodH $, consider it as an $H^*$-module, (\ref{eq10}) has the following form 
\begin{eqnarray}\Hom ^H (U(M), H) & \cong  &\Hom ^H (M, Rat(H^{*})) \nonumber\\
\nonumber & \cong &  \Hom ^H (M, (\lint ) \otimes H)\\
&\cong&\Hom_k(M,\lint),\label{eq101}\end{eqnarray} 
where the second isomorphism  follows from Sweedler's isomorphism (\ref{eq5}). 
Consequently, if $M$ has finite dimension
\begin{equation}\label{eq102} U(M)^*\cong \lint\ot M^*\end{equation}
(as vector spaces). Since $\lint\neq 0$, $U(M)\neq 0$ for any $M\neq 0$. If there exists
a comodule $M$ such that $\dim_kU(M)\leq \dim_kM$ then $\lint$ is one-dimensional.

\begin{lemma} \label{th5}
Let $H$ be a Hopf algebra over a field $k$ with a non-zero
 right integral then there exists a finite dimensional comodule $N$ such 
that $\dim U(N)\leq \dim_kN$.
\end{lemma}
\noindent
{\it Proof}. 
Let $C$ be a finite dimensional subcoalgebra of $H$. Then the action of $H^*$ on 
$\Hom_{H^*}(H,k)$ induces in the natural way an action on $\Hom_{H^*}(C,k)$:
$$\xi * g(c) = \sum_{(c)} \xi (c_{(1)}) g(c_{(2)}), \quad c\in C, g\in \Hom_{H^*}(C,k).$$
Moreover, there exists a left $H^*$-module morphism from $\Hom _{H^*} (H, k)$ 
to $\Hom _{H^*} (C, k)$  given by restriction:
$$\varphi _C : \Hom _{H^*} (H, k) \longrightarrow
 \Hom _{H^*} (C, k), f \mapsto \varphi _{C} (f)= f |_{C} .$$

We show that $\Hom_{H^*}(C,k)$ is rational.
Consider $C^*$ as a right $H$-comodule with the coaction given by condition 
$$\sum_{(g)} g_{(0)} (c) \otimes g_{(1)} = \sum_{(c)} c_{(1)} \otimes g(c_{(2)}),\quad
 \forall c \in C, g\in C^*.$$
Thus  $C^*$ is a rational left $H^*$-module, with the explicit action
$$(\xi * g)(c) = \sum_{(c)} \xi (c_{(1)}) g(c_{(2)}).$$
We see that the natural inclusion $\Hom_{H^*}(C,k)\longrightarrow
 \Hom_k(C,k)=C^*$ is compatible
with the actions of $H^*$, thus $\Hom_{H^*}(C,k)$ is an $H^*$-subcomodule of $C^*$,
whence rational.

Since $H$ is the union of its finite dimensional subcoalgebras, for
 a non-zero right integral $\chi$, there exists  a finite dimensional coalgebra $C$
such that $\chi|_C\not= 0$. Thus there exists an exact sequence 
\begin{equation}\label{eq_gamma0}
0 \longrightarrow K \longrightarrow \Hom _{H^*}
 (H, k) \longrightarrow \Gamma  \longrightarrow 0 \end{equation}
with $K=Ker \varphi _C , \Gamma =Im \varphi _C.$ 
$\Gamma$ is a rational $H^*$-module, being a submodule 
of $\Hom _{H^*} (C, k)$. $\Gamma  \neq 0$ since it 
contains $\varphi_C(\chi )$. Since $U$ is exact, the following sequence is also exact  
\begin{equation}\label{eq12}
0 \longrightarrow U(K) \longrightarrow U(\Hom _{H^*} (H, k))
 \longrightarrow U(\Gamma ) \longrightarrow 0. \end{equation}
By (\ref{eq7}),  $U (\Hom _{H^*} (H, k)) \cong k$, hence 
\begin{equation} \label{eq13}0 \longrightarrow U(K) \longrightarrow 
 k \longrightarrow U(\Gamma ) \longrightarrow 0.  \end{equation}
Thus, $\dim _k U(\Gamma ) =1\leq\dim_k\Gamma$. \eee

\begin{thm}\label{main_thm} Assume that the Hopf algebra $H$ possesses a non-zero
right integral then it is uniquely determined upto a constant. 
\end{thm}
\noindent{\it Proof.}
According to Lemma \ref{th5} and the discussion preceeding it, $\lint$ is one-dimensional. 
Using the symmetry between the notions of left and right integrals, we see that $\rint$ is
also one-dimensional.\eee
\begin{corollary}\label{lem_gamma}The space of right integrals $\rint\cong\Hom^H(H,k)\cong\Hom_{H^*}(H,k)$ with
the action of $H^*$ given in (\ref{eq8}) is a rational module.\end{corollary}
\noindent{\it Proof.} In the exact sequence (\ref{eq_gamma0}) $\Hom_{H^*}(H,k)$ 
is one-dimensional and $\Gamma\neq 0$, hence
is one-dimensional, too. Thus $K=0$ and $\Gamma\cong \Hom_{H^*}(H,k)$, whence
$\Hom_{H^*}(H,k)$ is rational.
\eee

Note that this one dimensional module induces a distinguished group-like element
introduced by Radford \cite{radford1}. Call this group-like element $\gamma$, we
have by definition, for any right integral $\chi$ on $H$, $\delta(\chi)=\chi\otimes \gamma$.
Since $\delta$ is induced from the action of $H^*$ on $\Hom_{H^*}(H,k)$, we have
the following equation for $\gamma$:
\begin{eqnarray*}
\xi*\chi(a)&=&\xi(\gamma)\chi(a)=\sum_{(a)}\xi(a_{(1)})\chi(a_{(2)}).
\end{eqnarray*}
or equivalently
\begin{equation}\label{eq_gamma}
\chi(a)\gamma=\sum_{(a)}a_{(1)}\chi(a_{(2)}).\end{equation}

\section{The bijectivity of the antipode}
In this section we prove that if $\rint\neq 0$ then the antipode is bijective.
\begin{lemma}\label{lem10}
Let H be a Hopf algebra, N be a  right $H$-comodule of finite dimension. Then 
$\Hom ^H (H, N)$ is a rational $H^*$-module and isomorphic to $N^{**}\otimes \Gamma$
as $H$-comodules.\end{lemma}
{\it Proof.} According to (\ref{eq15}), we have, for $N$ finite dimensional,
\begin{equation}\Hom^H (H, N) \cong N^{**} \otimes \Hom^H (H, k).
\end{equation}
It is more convenient to describe the inverse of this isomorphism, which is given by
$\phi : n \ot {f} \longmapsto  {f}_n , {f} _n (h) :=\sum_{(n)} n_{(0)} {f}(S(n_{(1)}) h). $

According to Corollary \ref{lem_gamma}, $N^{**} \otimes \Hom^H( H, k)$ 
is a right $H$-comodule, hence a left $H^*$-module 
with the action: 
$$\xi * (n \otimes {f}) := \sum_{(n)}\xi(S^2 (n_{(1)}) \gamma ))n_{(0)}  \ot  f   $$
for any $\xi $ in $H^*$ and  $n \ot {f} $ in $N^{**} \ot \Hom ^H (H, k).$ 

We prove that $\phi$ is a morphism of left $H^*$-modules (with respect to the above
actions). We have
\begin{eqnarray*}
\phi (\xi * (n \ot {f} )) (h)& = & \phi (n_{(0)} \ot f (\xi* (S^2 (n_{(1)})\Gamma ))(h) \\
&=&  \xi (S^2 (n_{(1)}) \gamma ) {f}_{n_{(0)}}(h) \\
&=& \sum_{(n)}n_{(0)} {f} (S(n_{(1)})h).  \xi (S^2 (n_{(2)})\gamma )
\end{eqnarray*} 
and 
\begin{eqnarray*}
\xi* (\phi( {f}\ot n)) (h) &=& \sum_{(h)(n)}\xi (h_{(1)}) n_{(0)} {f} (S(n_{(1)}) h_{(2)}). 
\end{eqnarray*}
According to (\ref{eq_gamma}), we have 
$f(S(n_{(1)})h ) \gamma  = \sum_{(n),(h)}S(n_{(2)})h_{(1)} {f} (S(n_{(1)})h_{(2)})$, whence
\begin{eqnarray*}
\sum_{(h)}h_{(1)} {f}(S(n_{(1)}) h_{(2)}) &=& \sum_{(n)}S^2 (n_{(2)})\gamma {f} (S(n_{(1)}) h). 
\end{eqnarray*}
Therefore, 
$$\sum_{(n)}n_{(0)} {f} (S(n_{(1)})h) \xi (S^2 (n_{(2)})) =
 \sum_{(h)(n)}\xi (h_{(1)}) n_{(0)} {f} (S(n_{(1)}) h_{(2)}).$$
Thus $\phi $ is a morphism of left $H^*$-modules, that is
 \begin{equation}  \Hom^H (H, N)  \cong N^{**} \otimes \Hom^H (H, k)
 \label {eq17}
\end{equation} as left $H^*$-modules.    
Since $N^{**} \otimes \Hom^H (H, k)$ is a rational $H^*$-module,  so is the module
$\Hom^H (H, N)$. \eee

\begin{pro} \label{pro7}
Let $H$ be a Hopf algebra with a non-zero integral
 and $U$ be the functor defined as in Proposition \ref{pro4}.
Then for any finite dimensional comodule $N$, we have
$U(N)^{**}\cong N\ot \Gamma$.
\end{pro}
{\it Proof.}
 From (\ref{eq9}) we have, for finite dimensional right $H$-comodules $ M, N$: 
$$\Hom^H (U(N), M) \cong \Hom^H (N, \Hom^H (H, M))$$
From (\ref{eq17}),
$$ \Hom^H (N, \Hom^H (H, M)) \cong \Hom^H (N, M^{**} \otimes   \Gamma ).$$
We have $\Gamma $ is invertible ($\Gamma^*\ot\Gamma\cong k\cong\Gamma^*\ot\Gamma$), hence 
$\Hom^H (M \otimes \Gamma, N \otimes \Gamma ) \cong \Hom^H (M, N)$. Therefore 
$$\Hom^H (N, M^{**} \otimes \Gamma ) \cong \Hom^H (N \otimes \Gamma ^*,
 M^{**} ).$$
Using (\ref{eq14}) and (\ref{eq15}), we have, for any finite dimensional comodule $M$,
\begin{eqnarray*}\Hom^H (U(N) ^{**} \otimes M^* , k)& 
\cong&\Hom^H(U(N),M)\\
&\cong& \Hom^H(N,M^{**}\ot \Gamma)\\
&\cong &\Hom^H(N\ot\Gamma^*,M^{**})\\
&\cong &  \Hom^H (N \otimes \Gamma \otimes M^* , k)\end{eqnarray*}
hence $ U(N)^{**} \cong N \otimes \Gamma $. 
 \eee 

\begin{corollary}
 Let $H$ be a Hopf algebra over a field $k$ with antipode 
$S$, and assume $\rint\ne 0$. Then $S$ is bijective.
\end{corollary}
{\it Proof.} We have $(\Gamma \otimes U(N)^*) ^{*}\cong U(N)^{**}\ot\Gamma^*\cong N$, i.e., every 
 finite dimensional right $H$-comodule $N$ has right dual.  Hence the antipode
is bijective. \eee
\begin{corollary} Let $H$ possess a non-zero integral. Then the injective 
envelope of $k$ is finite dimensional.\end{corollary}
\noindent{\it Proof.}  Let $\phi$ be a non-zero right integral on $H$.
According to the proof of Theorem \ref{main_thm}, there exists a finite
dimensional right ideal $J$ of $H$ such that the restriction of $\phi$ on it is non-zero.
Choose such a right ideal $J$
with minimal dimension. We show that $J$ is a direct summand of $H$ as $H$-right
comodule. Let $\phi_J$ be the restriction on $J$: $\phi_J:J\longrightarrow k$. 

According to Propostion \ref{pro7}, the functor $\Hom^H(H,-)$ is exact
on finite dimensional comodules, hence there exists a morphism $\pi:H\longrightarrow
J$, such that $\phi=\phi_J\circ \pi$. Hence $\phi_J=\phi_J\circ\pi|_J.$ By the minimality
of $J$, $\pi$ should be surjective, whence $\pi|_J$ should be bijective since $J$ is
finite dimensional. Therefore $J$ is a direct summand of $H$, whence $J$ is projective
with respect to finite dimensional comodules. Since each comodule is the union of its
finite dimensional subcomodules, we conclude that $J$ is projective in $\ComodH$. 
Hence $J^*$ is injective, and thus, the injective envelope of $k$. \eee

\end{document}